\documentstyle{amltd2004}
\begin{document}
\annalsline{158}{2003}
\received{January 24, 2002}
\startingpage{323}
\def\bye{\end{document}}
 \font\tenrm=cmr10
\def\ritem#1{\item[{\rm #1}]}
\input amssym.def
\input amssym.tex
\def\hensp#1{\quad\hbox{#1}\quad }
\catcode`\@=11
\font\twelvemsb=msbm10 scaled 1100
\font\tenmsb=msbm10
\font\ninemsb=msbm10 scaled 800
\newfam\msbfam
\textfont\msbfam=\twelvemsb  \scriptfont\msbfam=\ninemsb
  \scriptscriptfont\msbfam=\ninemsb
\def\msb@{\hexnumber@\msbfam}
\def\Bbb{\relax\ifmmode\let\next\Bbb@\else
 \def\next{\errmessage{Use \string\Bbb\space only in math
mode}}\fi\next}
\def\Bbb@#1{{\Bbb@@{#1}}}
\def\Bbb@@#1{\fam\msbfam#1}
\catcode`\@=12

 \catcode`\@=11
\font\twelveeuf=eufm10 scaled 1100
\font\teneuf=eufm10
\font\nineeuf=eufm7 scaled 1100
\newfam\euffam
\textfont\euffam=\twelveeuf  \scriptfont\euffam=\teneuf
  \scriptscriptfont\euffam=\nineeuf
\def\euf@{\hexnumber@\euffam}
\def\frak{\relax\ifmmode\let\next\frak@\else
 \def\next{\errmessage{Use \string\frak\space only in math
mode}}\fi\next}
\def\frak@#1{{\frak@@{#1}}}
\def\frak@@#1{\fam\euffam#1}
\catcode`\@=12

\title{Global existence and convergence for a\\ higher order flow in conformal 
geometry}
\shorttitle{Global existence} 
  \author{Simon Brendle}
 \institutions{Princeton University, Princeton, NJ\\
{\eightpoint {\it  E-mail address\/}:  brendle@math.princeton.edu
}}

\newcommand{\tr}{\hbox{\rm tr}\,}

\section{Introduction}

An important problem in conformal geometry is the construction of 
conformal metrics for which a certain curvature quantity equals a prescribed function, e.g. a
constant. In two dimensions, the uniformization theorem assures the
existence of a conformal metric with constant Gauss curvature. Moreover, 
J.\ Moser \cite{Mo} proved that for every positive function $f$ on $S^2$
satisfying $f(x) = f(-x)$ for all $x \in S^2$ there exists a conformal metric on
$S^2$ whose Gauss curvature is equal to $f$. 

A natural 
conformal invariant in dimension four is $$ Q = 
-\frac{1}{6} \, (\Delta R - R^2 + 3 \, |{\rm Ric}|^2),$$  
where $R$ denotes the scalar curvature and ${\rm Ric}$ the Ricci tensor. 
This formula can also be
written in the form $$ Q = -\frac{1}{6} \, (\Delta R - 6 \, \sigma_2(A)),$$  
where $$ A = {\rm Ric} - \frac{1}{6} \, Rg$$  is the Schouten tensor of $M$ and 
$$ \sigma_2(A) = \frac{1}{2} \, (\tr A)^2 - \frac{1}{2} \, |A|^2$$  is the
second elementary symmetric polynomial in its eigenvalues. 
Under a conformal change of the metric $$ g = e^{2w} g_0,$$  the quantity $Q$
transforms according to $$ Q = e^{-4w} (Q_0 + P_0 w),$$  where $P_0$ denotes the
Paneitz operator with respect to $g_0$. The Gauss-Bonnet-Chern theorem asserts that 
$$ \int_M Q \, dV + \int_M \frac{1}{4} \, |W|^2 \, dV = 8\pi^2 \chi(M).$$  
Since the Weyl tensor $W$ is conformally invariant, it follows 
that the expression $$ \int_M Q \, dV$$  is conformally invariant, too. 
The quantity $Q$ plays an important role in four-dimensional conformal geometry;
see \cite{BCY}, \cite{CY1}, \cite{CGY}, \cite{Gu1} (note that our notation differs slightly from that 
in \cite{BCY}, \cite{CY1}). Moreover, the Paneitz operator plays a similar role as the 
Laplace operator in dimension two; compare \cite{BCY}, \cite{CY1},
\cite{CGY}, \cite{DMA1}, \cite{DMA2}. We also 
note that the Paneitz operator is of considerable interest in mathematical 
physics, see \cite[\SS IV.4]{Co}.  

T.\ Branson, S.-Y.\ A.\ Chang and P.\ Yang \cite{BCY} studied metrics for which the
curvature quantity $Q$ is constant. Since $$ \int_M Q \, dV$$  is conformally
invariant, these metrics minimize the functional $$ \int_M Q^2 \, dV$$  among all
conformal metrics with fixed volume. In addition, these metrics are critical points of the 
functional $$ E_1[w] = \int_M 2 \, w \, P_0 w \, dV_0 + \int_M 4 \, Q_0 \, w \, dV_0 - 
\int_M Q_0 \, dV_0 \,
\log \bigg ( \int_M e^{4w} \, dV_0 \bigg ),$$  where $g_0$ denotes a fixed 
metric on $M$ and $g = e^{2w} g_0$. 

According to the results in [2], one can construct conformal metrics of constant $Q$-curvature by minimizing the
functional $E_1[w]$ provided that the Paneitz operator is weakly positive and the integral of the $Q$-curvature on
$M$ is less than that on the standard sphere $S^n$.  In dimension four, M. Gursky \cite{Gu2} proved that both conditions are
satisfied if 
$$
Y(g_0)\ge 0
$$
and $$\int_M Q_0\, dV_0\ge 0,
$$
and $M$ is not conformally equivalent to the standard sphere $S^4$.

C.\ Fefferman and R.\ Graham \cite{FG1}, \cite{FG2} 
established the existence of a conformally invariant self-adjoint 
operator with leading term $(-\Delta)^{\frac{n}{2}}$ 
in all even dimensions $n$. Moreover, there is a
curvature quantity which transforms according to 
$$ Q = e^{-nw} (Q_0 + P_0 w)$$  for $$ g = e^{2w} g_0.$$  This implies that 
the expression $$ \int_M Q \, dV$$  is 
conformally invariant. Hence, a metric with $Q = \hbox{\rm constant}$ minimizes 
the functional $$ \int_M Q^2 \, dV$$  among all
conformal metrics with fixed volume. Finally, the analogue of the functional
$E_1[w]$ is given by 
$$ E_1[w] = \int_M \frac{n}{2} \, w \, P_0 w \, dV_0 + \int_M n \, Q_0 \, w \, dV_0 - 
\int_M Q_0 \, dV_0 \,
\log \bigg ( \int_M e^{nw} \, dV_0 \bigg ).$$

Our aim is to construct conformal metrics for which 
the curvature quantity $Q$ is a constant multiple of a prescribed positive
function $f$ on $M$. This equation is the Euler-Lagrange equation for 
the functional 
$$ E_f[w] = \int_M \frac{n}{2} \, w \, P_0 w \, dV_0 + \int_M n \, Q_0 \, w \, dV_0 - 
\int_M Q_0 \, dV_0 \,
\log \bigg ( \int_M e^{nw} \, f \, dV_0 \bigg ).$$  
We construct critical points of the functional $E_f[w]$ using 
the gradient flow for $E_f[w]$. A similar method was used by 
R.\ Ye \cite{Ye} to prove Yamabe's theorem for locally conformally flat 
manifolds. K.\ Ecker and G.\ Huisken \cite{EH} 
used a variant of mean curvature flow to construct hypersurfaces with 
prescribed mean curvature in cosmological spacetimes. 

The flow of steepest descent for the functional $E_f[w]$ is given by 
$$ \frac{\partial}{\partial t} g = 
- \Big ( Q - \frac{\overline{Q} \, f}{\overline{f}} \Big ) \, g.$$  Here, 
$\overline{Q}$ and $\overline{f}$ denote the mean values of $Q$ and $f$ respectively, 
i.e.\ $$ \int_M (Q - \overline{Q}) \, dV = 0\quad\hbox{ and }\quad  \int_M (f - \overline{f}) \, dV = 0.$$  
This evolution equation preserves the conformal structure of $M$. Moreover, 
since $$ \int_M \Big ( Q - \frac{\overline{Q} \, f}{\overline{f}} \Big ) \, dV 
= \int_M \Big ( \overline{Q} - \frac{\overline{Q} \, \overline{f}}{\overline{f}} \Big ) \, dV = 0,$$  
the volume of $M$ 
remains constant. From this it follows that $\overline{Q}$ is constant in 
time. If we write $g = e^{2w} g_0$ for a fixed metric $g_0$, then
the evolution equation takes the form 
$$ \frac{\partial}{\partial t} w = -\frac{1}{2} \, e^{-nw} \, P_0 w 
- \frac{1}{2} \, e^{-nw} \, Q_0 + \frac{1}{2} \, \frac{\overline{Q} \, f}{\overline{f}}, \pagebreak $$  
where $P_0$ denotes the Paneitz operator with respect to $g_0$. 
Therefore, the function $w$ satisfies a 
quasilinear parabolic equation of order $n$ involving the critical Sobolev exponent. 
Moreover, the reaction term is nonlocal, since $\overline{f}$ involves values of $w$ on the 
whole of $M$.

\proclaim{Theorem}\label{convergence.1}
Assume that the Paneitz operator $P_0$ is weakly positive with kernel consisting 
of the constant functions. Moreover{\rm ,} assume that 
$$ \int_M Q_0 \, dV_0 < (n-1)! \, \omega_n.$$  Then the evolution  equation $$ \frac{\partial}{\partial t} g 
= - \Big ( Q - \frac{\overline{Q} \, f}{\overline{f}} \Big ) \, g$$  has a 
solution which is defined for all times and converges to a metric with $$ \frac{Q}{f} = 
\frac{\overline{Q}}{\overline{f}}.$$ 
\endproclaim

On the standard sphere $S^n$, we have 
$$ \int_M Q \, dV = (n-1)! \, \omega_n\, ;$$  hence Theorem \ref{convergence.1} 
cannot be applied. In fact, the conclusion of Theorem \ref{convergence.1} 
fails for $M = S^n$. To see this, one can consider the Kazdan-Warner identity 
$$ \int_{S^n} \langle \nabla_0 Q,\nabla_0 x_j \rangle \, e^{nw} \, dV_0 = 0\, ;$$  
see \cite{CY1}. If
$f$ is an increasing function of $x_j$, then $$ \int_{S^n} \langle \nabla_0 Q,
\nabla_0 x_j \rangle \, e^{nw} \, dV_0 > 0.$$  Consequently, there is no conformal
metric on $S^n$ satisfying $$ \frac{Q}{f} = 
\frac{\overline{Q}}{\overline{f}}.$$  
Nevertheless, the conclusion of Theorem \ref{convergence.1} holds if 
$f(x) = f(-x)$ and $w(x) = w(-x)$ for all 
$x \in S^n$. This is a generalization of Moser's theorem \cite{Mo}.

\proclaim{Theorem}\label{convergence.2}
Suppose that $M = {\bf RP}^n$. Then the evolution equation $$ \frac{\partial}{\partial t} g 
= - \Big ( Q - \frac{\overline{Q} \, f}{\overline{f}} \Big ) \, g$$  has a 
solution which is defined for all times and converges to a metric with $$ \frac{Q}{f} = 
\frac{\overline{Q}}{\overline{f}}.$$ 
\endproclaim

Combining Theorem 1.2 with M. Gursky's result \cite{Gu2} gives 

\proclaim{Theorem}\label{convergence.3}
Suppose that $M$ is a compact manifold of dimension four satisfying 
$$ Y(g_0) \geq 0\hensp{and} \int_M Q_0 \, dV_0 \geq 0.$$  
Moreover{\rm ,} assume that $M$ is not conformally equivalent to the standard 
sphere~$S^4$. Then the evolution equation $$ \frac{\partial}{\partial t} g 
= - \Big ( Q - \frac{\overline{Q} \, f}{\overline{f}} \Big ) \, g$$  has a 
solution which is defined for all times and converges to a metric with $$ \frac{Q}{f} = 
\frac{\overline{Q}}{\overline{f}}.$$ 
\endproclaim

Finally, we prove a compactness theorem for conformal metrics on $S^n$. 
In two dimensions, the corresponding result was first established by X. Chen 
\cite{Ch1} (see also \cite{St}). 

\proclaim{Proposition}\label{compactness.result}
Let $g_k = e^{2w_k} \, g_0$ be a sequence of conformal metrics on $S^n$ with fixed
volume such that $$ \int_{S^n} Q_k^2 \, dV_k \leq C.$$  
Assume that for every point $x \in S^n$ there 
exists $r > 0$ such that 
$$ \lim_{r \to 0} \lim_{k \to \infty} \int_{B_r(x)} |Q_k| \, dV_k < \frac{1}{2}
\, (n-1)! \, \omega_n.$$  
Then the sequence $w_k$ is uniformly bounded in $H^n$.
\endproclaim

The evolution equation can be viewed as a generalization of the Ricci flow on 
compact surfaces. In dimension four, the quantity $Q$ plays a similar role as the 
Gauss curvature in dimension two. Moreover, the energy 
functional $E_1[w]$ corresponds to 
the Liouville energy studied by B.\ Osgood, R.\ Phillips and 
P.\ Sarnak in \cite{OPS}. 

It was shown by R.\ Hamilton \cite{Ha} and B.\ Chow \cite{Ch2} that every 
solution of the Ricci flow on a compact surface exists for all time 
and converges exponentially to a
metric with constant Gauss curvature. A different approach 
was introduced by X.\ Chen \cite{Ch1} in his work on the Calabi flow. Similar methods were used by M.\ Struwe
\cite{St} to establish global existence  and exponential convergence for the Ricci flow on compact surfaces, and by X.
Chen and G. Tian [7] to prove convergence of the K\"ahler-Ricci flow on K\"ahler-Einstein surfaces. For the Ricci flow,
the situation is more complicated since the Calabi energy is not decreasing along the flow. H.\ Schwetlick \cite{Sc}
used similar arguments to  deduce global existence and convergence for a natural 
sixth order flow on surfaces. The approach used \pagebreak in \cite{Ch1} and \cite{St} 
is based on integral estimates and does not rely on the maximum principle. These 
ideas are also useful in our situation. This is due to the fact that the 
equation studied in this paper has higher order, hence the maximum principle is not available.  

In Section 2 we derive the evolution equation for the conformal factor and the
curvature quantity $Q$. In
Section 3 we show that the solution is bounded in $H^{\frac{n}{2}}$. In Sections 4 and 5 we 
show that the solution exists for all time, and in Section 6 we prove that the
evolution equation converges to a stationary solution. Finally, the proof of 
Proposition \ref{compactness.result} is carried out in Section 8. 

The author would like to thank S.-Y.\ A.\  Chang and J.\ Viaclovsky for 
helpful comments.

\section{The evolution equations for $w$ and $Q - \frac{\overline{Q} \,
f}{\overline{f}}$}

Since the evolution equation preserves the conformal structure, we may write $g = e^{2w} \, g_0$
for a fixed metric $g_0$ and some real-valued function $w$. Then we have the 
formula $$ Q = e^{-nw} \, (Q_0 + 
P_0w),$$  where $P_0$ denotes the Paneitz operator with respect 
to the metric $g_0$. Hence, the function $w$ obeys the evolution equation 
$$ \frac{\partial}{\partial t} w = -\frac{1}{2} \, e^{-nw} \, P_0 w 
- \frac{1}{2} \, e^{-nw} \, Q_0 + \frac{1}{2} \, 
\frac{\overline{Q} \, f}{\overline{f}}.$$  
Differentiating both sides with respect to $t$, we obtain 
$$ \frac{\partial}{\partial t} \Big ( Q - \frac{\overline{Q} \, f}{\overline
{f}} \Big ) = -\frac{1}{2} \, P \Big ( Q 
- \frac{\overline{Q} \, f}{\overline{f}} \Big ) + \frac{n}{2} \, Q \Big ( 
Q - \frac{\overline{Q} \, f}{\overline{f}} \Big ) + \frac{\overline{Q} \,
f}{\overline{f} \, ^2} \, \frac{d}{dt} \overline{f},$$  where $P = e^{-nw} \, P_0$ is the Paneitz 
operator with respect to the metric $g$. It follows from the evolution
equation for $w$ that 
$$ \frac{d}{dt} \overline{f} = -\int_M \frac{n}{2} \, f \, \Big ( Q - \frac
{\overline{Q} \, f}{\overline{f}} \Big ) \, dV.$$  
This implies 
\begin{eqnarray*} 
\frac{\partial}{\partial t} \Big ( Q - \frac{\overline{Q} \, f}{\overline
{f}} \Big ) &=& -\frac{1}{2} \, P \Big ( Q 
- \frac{\overline{Q} \, f}{\overline{f}} \Big ) + \frac{n}{2} \, Q \Big ( 
Q - \frac{\overline{Q} \, f}{\overline{f}} \Big ) 
\\[5pt] &&- \frac{n}{2} \, \frac{\overline{Q} \, f}{\overline{f}} \int_M \frac{f}{\overline{f}} \, 
\Big ( Q - \frac{\overline{Q} \, f}{\overline{f}} \Big ) \, 
dV, \end{eqnarray*} 
where $P$ denotes the Paneitz operator with respect to 
the \pagebreak metric $g$.

\section{Boundedness of $w$ in $H^{\frac{n}{2}}$}

We consider the functional 
$$ E_f[w] = \int_M \frac{n}{2} \, w \, P_0 w \, dV_0 + \int_M n \, Q_0 \, w \, dV_0 - 
\int_M Q_0 \, dV_0 \,
\log \bigg ( \int_M f \, e^{nw} \, dV_0 \bigg ).$$  
Since $P_0$ is self-adjoint,  
\begin{eqnarray*} 
\frac{d}{dt} E_f[w] 
&=& \int_M \frac{n}{2} \, \frac{\partial}{\partial t} w \, P_0 w \, dV_0 
+ \int_M \frac{n}{2} \, w \, P_0 \frac{\partial}{\partial t} w \, dV_0 
+ \int_M n \, Q_0 \, \frac{\partial}{\partial t} w \, dV_0 
\\ &&- \int_M n \, \frac{\overline{Q} \, f}{\overline{f}} \, \frac{\partial}{\partial
t} w \, dV \\ 
&=& \int_M n \, P_0 w \, \frac{\partial}{\partial t} w \, dV_0 
+ \int_M n \, Q_0 \, \frac{\partial}{\partial t} w \, dV_0 
 - \int_M n \, \frac{\overline{Q} \, f}{\overline{f}} \, \frac{\partial}{\partial
t} w \, dV \\ 
&=& \int_M n \, Q \, \frac{\partial}{\partial t} w \, dV 
 - \int_M n \, \frac{\overline{Q} \, f}{\overline{f}} \, \frac{\partial}{\partial
t} w \, dV \\ 
&=& \int_M n \, \Big ( Q - \frac{\overline{Q} \, f}{\overline
{f}} \Big ) \, \frac{\partial}{\partial t} w \,
dV. 
\end{eqnarray*}  
Since the time derivative of $w$ is given by $$ \frac{\partial}{\partial t} w =
-\frac{1}{2} \, \Big ( Q - \frac{\overline{Q} \, f}{\overline{f}} \Big ),$$  we obtain 
$$ \frac{d}{dt} E_f[w] = -\int_M \frac{n}{2} \, \Big ( Q - \frac{\overline{Q} 
\, f}{\overline{f}} \Big )^2 \, dV.$$  
In particular, the functional $E_f[w]$ is decreasing under the evolution
equation. This implies 
$$ E_f[w] \leq C.$$  
In the first step, we consider the case $$ \int_M Q_0 \, dV_0 < 0.$$  
Using Jensen's inequality we obtain 
$$ \log \bigg ( \int_M e^{n(w - \overline{w})} \, dV_0 \bigg ) \geq -C.$$  
This implies 
\begin{eqnarray*} 
E_f[w] &\geq& \int_M \frac{n}{2} \, w \, P_0 w \, dV_0 + \int_M n \, Q_0 \, w \, dV_0 \\ 
&&- \int_M Q_0 \, dV_0 \,
\log \bigg ( \int_M e^{nw} \, dV_0 \bigg ) - C \end{eqnarray*}
\begin{eqnarray*}
&=& \int_M \frac{n}{2} \, w \, P_0 w \, dV_0 + \int_M n 
\, Q_0 \, (w - \overline{w}) \, dV_0 \\ 
&&- \int_M Q_0 \, dV_0 \,
\log \bigg ( \int_M e^{n(w-\overline{w})} \, dV_0 \bigg ) - C \\ 
&\geq&\int_M \frac{n}{2} \, w \, P_0 w \, dV_0 + \int_M n \, Q_0 \, (w - \overline{w}) \, dV_0
- C \\ 
&\geq &2\delta \int_M \big ( (-\Delta_0)^{\frac{n}{4}} w \big )^2 \, dV_0 + \int_M n \, Q_0 \, (w -
\overline{w}) \, dV_0 - C \\ 
&\geq& \delta \int_M \big ( (-\Delta_0)^{\frac{n}{4}} w \big )^2 \, dV_0 - C. 
\end{eqnarray*} 
In the second step, we consider the case 
$$ 0 \leq \int_M Q_0 \, dV_0 < (n - 1)! \, \omega_n.$$  
Since the Paneitz operator $P_0$ is self-adjoint and weakly positive, 
it has a square root $P_0^{\frac{1}{2}}$. Moreover, the kernel of
$P_0^{\frac{1}{2}}$ coincides with the kernel of $P_0$, which 
consists of the constant functions. Thus, we conclude that 
$$ w(y) - \overline{w} = \int_M P_0^{\frac{1}{2}} w(z) \, H(y,z) \, dV_0(z)$$  
for a suitable function $H(y,z)$. The leading term in the asymptotic expansion of 
the kernel $H(y,z)$ coincides with that of the Green's function for the operator
$(-\Delta)^{\frac{n}{4}}$ in ${\bf R}^n$. Hence, we can apply 
an inequality of D. Adams (see \cite[Theorems~1 and 2]{Ad}). This 
implies 
$$ \int_M e^{\frac{2^n \pi^n n}{\omega_{n-1}} \, \frac
{(w - \overline{w})^2}{\int_M (P_0^{\frac{1}{2}} w)^2 \, dV_0}} \, dV_0 \leq C,$$  
hence 
$$ \int_M e^{\frac{2^n \pi^n n}{\omega_{n-1}} \, \frac
{(w - \overline{w})^2}{\int_M w \, P_0 w \, dV_0}} \, dV_0 \leq C.$$  
Since $$ \omega_{n-1} \omega_n = \frac{2^{n+1} \pi^n}{(n - 1)!},$$  
we obtain $$ \int_M e^{n(w - \overline{w})} \, dV_0 
\leq C \, e^{\int_M \frac{n}{2(n-1)! \, \omega_n} \, w \, P_0 w \, dV_0}.$$  
From this it follows that
\begin{eqnarray*} 
E_f[w] &\geq& \int_M \frac{n}{2} \, 
w \, P_0 w \, dV_0 + \int_M n \, Q_0 \, w \, dV_0 \\ 
&&- \int_M Q_0 \, dV_0 \,
\log \bigg ( \int_M e^{nw} \, dV_0 \bigg ) - C \end{eqnarray*}
\begin{eqnarray*}
&=& \int_M \frac{n}{2} \, w \, P_0 w \, dV_0 + \int_M n \, Q_0 \, (w - \overline{w}) \, dV_0 \\ 
&&- \int_M Q_0 \, dV_0 \,
\log \bigg ( \int_M e^{n(w-\overline{w})} \, dV_0 \bigg ) - C \\ 
&\geq& \bigg ( 1 - \frac{\int_M Q_0 \, dV_0}{(n-1)! \, \omega_n} \bigg ) \, \int_M \frac
{n}{2} \, w \, P_0 w
\, dV_0 + \int_M n \, Q_0 \, (w - \overline{w}) \, dV_0 - C \\ 
&\geq &2\delta \int_M \big ( (-\Delta_0)^{\frac{n}{4}} w \big )^2 \, dV_0 + \int_M 4 \, Q_0 \, (w -
\overline{w}) \, dV_0 - C \\ 
&\geq &\delta \int_M \big ( (-\Delta_0)^{\frac{n}{4}} w \big )^2 \, dV_0 - C. 
\end{eqnarray*} 
Since $E_f[w]$ is bounded from above, we conclude that 
$$ \int_M \big ( (-\Delta_0)^{\frac{n}{4}} w \big )^2 \, dV_0 \leq C;$$  
hence 
$$ \|w - \overline{w}\|_{H^{\frac{n}{2}}} \leq C.$$  
Using an inequality of N. Trudinger, we obtain 
$$ \int_M e^{\alpha(w - \overline{w})} \, dV_0 \leq C$$  
for all real numbers $\alpha$. 
In particular, we have 
$$ \int_M e^{n(w - \overline{w})} \, dV_0 \leq C.$$  
Since $  \int_M e^{nw} \, dV_0 = 1,$   we conclude that 
$  e^{-n\overline{w}} \leq C;$   hence 
 $ -C \leq \overline{w} \leq C.$   
This implies 
 $ \|w\|_{H^{\frac{n}{2}}} \leq C$ 
and 
 $ \int_M e^{\alpha w} \, dV_0 \leq C$  
for all real  numbers $\alpha$. 
Since the functional $E_f[w]$ is bounded from below, we finally obtain 
$$ \int_0^T \int_M \Big ( Q - \frac{\overline{Q} \, f}{\overline{f}} \Big )^2 \, dV \, dt \leq C.$$

\section{Boundedness of $w$ in $H^n$ for $0 \leq t \leq T$}

Let $T$ be a fixed, positive real number. We claim that 
$$ \|w\|_{H^n} \leq C$$  
for all $0 \leq t \leq T$. For the sake of brevity, we put 
\begin{eqnarray*} v 
&=& -\frac{1}{2} \, e^{\frac{nw}{2}} \, \Big ( Q - \frac{\overline{Q} 
\, f}{\overline{f}} \Big )  = 
e^{\frac{nw}{2}} \, \frac{\partial}{\partial t} w \\ &=& -\frac{1}{2} \, 
e^{-\frac{nw}{2}} \, P_0 w 
- \frac{1}{2} \, e^{-\frac{nw}{2}} \, Q_0 
+ \frac{1}{2} \, e^{\frac{nw}{2}} \, \frac{\overline{Q} \,
f}{\overline{f}}. \end{eqnarray*} 
This implies 
 $$ \frac{\partial}{\partial t} w = e^{-\frac{nw}{2}} \, v\hensp{   and  }
  P_0 w = -2 \, e^{\frac{nw}{2}} \, v - Q_0 + e^{nw} \, \frac{\overline{Q} \, f}{\overline{f}}.$$ 
From this it follows that 
\begin{eqnarray*} 
\frac{d}{dt} \bigg ( \int_M (P_0 w)^2 \, dV_0 \bigg ) 
&= &-\int_M 4 \, (e^{\frac{nw}{2}} \, v) \, P_0(e^{-\frac{nw}{2}} \, v) \, dV_0 \\ 
&&- \int_M 2 \, Q_0 \, P_0(e^{-\frac{nw}{2}} \, v) \\ 
&&+ \int_M \frac{2 \, \overline{Q}}{\overline{f}} \, (e^{nw} \, f) \, P_0(e
^{-\frac{nw}{2}} \, v) \, dV_0. \end{eqnarray*} 
This implies 
\begin{eqnarray*} 
\frac{d}{dt} \bigg ( \int_M (P_0 w)^2 \, dV_0 \bigg ) 
&=& -\int_M 4 \, (-\Delta_0)^{\frac{n}{4}}(e^{\frac{nw}{2}} \, v) \, 
(-\Delta_0)^{\frac{n}{4}}(e^{-\frac{nw}{2}} \, v) \, dV_0 \\ 
&&- \int_M 2 \, P_0 Q_0 \, (e^{-\frac{nw}{2}} \, v) \\ 
&&+ \int_M \frac{2 \, \overline{Q}}{\overline{f}} \, (-\Delta_0)^{\frac{n}{4}} 
(e^{nw} \, f) \, (-\Delta_0)^{\frac{n}{4}}
(e^{-\frac{nw}{2}} \, v) \, dV_0 \\ &&+ \hbox{ \rm lower order terms}. \end{eqnarray*} 
Here, we adopt the convention that 
$$ (-\Delta_0)^{m+\frac{1}{2}} = \nabla_0 \, (-\Delta_0)^m$$  
for all integers $m$ (see \cite{Ad}). The right-hand side involves derivatives
of $v$ and $w$ of order at most $\frac{n}{2}$. Moreover, the total number of
derivatives is at most $n$. Therefore, we obtain 
\begin{eqnarray*} 
\frac{d}{dt} \bigg ( \int_M (P_0 w)^2 \, dV_0 \bigg ) 
&=& -\int_M 4 \, \big ( (-\Delta_0)^{\frac{n}{4}} v \big )^2 \, dV_0 \\ &&+ C 
\sum_{k_1, \ldots, k_m} \int_M |\nabla_0^{k_1} v| \cdot 
|\nabla_0^{k_2} v| \cdot |\nabla_0^{k_3} w| \cdots |\nabla_0^{k_m} w| \, dV_0 \\ &&+ C 
\sum_{l_1, \ldots, l_m} \int_M |\nabla_0^{l_1} v| \cdot 
|\nabla_0^{l_2} w| \cdots |\nabla_0^{l_m} w| \, e^{\alpha w} \, dV_0. 
\end{eqnarray*}
The first sum is taken over all $m$-tuples $k_1, \ldots, k_m$ with $m \geq 3$
satisfying the conditions 
\begin{eqnarray*}
0\le k_i\le \frac{n}{2} &&  \hbox{for } 1\le i\le 2,\\
1\le k_i\le \frac{n}{2}&&\hbox{for }3\le i\le m,
\end{eqnarray*}
and
$$
k_1+\cdots + k_m \le n.
$$
To estimate this term, we choose real numbers $p_1, \ldots, p_m \in [2,\infty[$ such that 
\begin{eqnarray*}
 k_i  \leq  \frac{n}{p_i} &&\hbox{  for } \quad 1 \leq i \leq 2,\\
 \frac{n}{p_i}  <  k_i &&\hbox{ for } \quad 3
\leq i \leq m
\end{eqnarray*}  and $$ \frac{1}{p_1} +\cdots +\frac{1}{p_m} = 1.$$  
Moreover, we define real numbers $\theta_1, \ldots, \theta_m$ by 
$$ \theta_i = \frac{k_i - \frac{n}{p_i} + \frac{n}{2}}{\frac{n}{2}} \in [0,1]\quad \hbox{
for }1 \leq i \leq 2$$ and 
$$ \theta_i = \frac{k_i - \frac{n}{p_i}}{\frac{n}{2}} \in \hspace{1.2mm} ]0,1[ \quad \hbox{ for }\quad 3 \leq i \leq
m.$$  Then we have 
$  \theta_1 +\cdots +\theta_m \leq 2$;
hence 
$$ \theta_3 +\cdots +\theta_m \leq (1-\theta_1) + (1-\theta_2).$$  
Since $w$ is bounded in $H^{\frac{n}{2}}$, this implies 
\begin{eqnarray*} 
&&-\int_M 2 \, \big ( (-\Delta_0)^{\frac{n}{4}} v \big )^2 \, dV_0 
\\ &&+ C \sum_{k_1, \ldots, k_m} \int_M |\nabla_0^{k_1} v| \cdot 
|\nabla_0^{k_2} v| \cdot |\nabla_0^{k_3} w| \cdots |\nabla_0^{k_m} w| \, dV_0 
\\ &\leq &-\|v\|_{H^{\frac{n}{2}}}^2 + C 
\sum_{k_1, \ldots, k_m} \|\nabla_0^{k_1} v\|_{L^{p_1}} \cdot 
\|\nabla_0^{k_2} v\|_{L^{p_2}} \cdot \|\nabla_0^{k_3} w\|_{L^{p_3}} 
\cdots \|\nabla_0^{k_m} w\|_{L^{p_m}} \\ 
&\leq& -\|v\|_{H^{\frac{n}{2}}}^2 + C 
\sum_{k_1, \ldots, k_m} \|v\|_{H^{k_1 - \frac{n}{p_1} + \frac{n}{2}}} \\&&\cdot 
\|v\|_{H^{k_2 - \frac{n}{p_2} + \frac{n}{2}}} \cdot 
\|w\|_{H^{k_3 - \frac{n}{p_3} + \frac{n}{2}}} 
\cdots \|w\|_{H^{k_m - \frac{n}{p_m} + \frac{n}{2}}} \\ 
&\leq &-\|v\|_{H^{\frac{n}{2}}}^2 + C 
\sum_{k_1, \ldots, k_m} \|v\|_{L^2}^{(1-\theta_1)+(1-\theta_2)} \, 
\\
&&\cdot \|v\|_{H^{\frac{n}{2}}}^{\theta_1+\theta_2} 
\, \|w\|_{H^{\frac{n}{2}}}^{(1-\theta_3)+\cdots +(1-\theta_m)} \, 
\|w\|_{H^n}^{\theta_3+\cdots +\theta_m} \\ 
&\leq& -\|v\|_{H^{\frac{n}{2}}}^2 + C 
\sum_{k_1, \ldots, k_m} \|v\|_{L^2}^{(1-\theta_1)+(1-\theta_2)} \, 
\|v\|_{H^{\frac{n}{2}}}^{\theta_1+\theta_2} 
\, \|w\|_{H^n}^{\theta_3+\cdots +\theta_m} \\ 
&\leq& C 
\sum_{k_1, \ldots, k_m} \|v\|_{L^2}^2 \, 
\|w\|_{H^n}^{\frac{2(\theta_3+\cdots +\theta_m)}{(1-\theta_1)+(1-\theta_2)}} \\ 
&\leq& C \, \|v\|_{L^2}^2 \, (\|w\|_{H^n}^2 + 1). 
\end{eqnarray*}  
The second sum is taken over all $m$-tuples $l_1, \ldots, l_m$ with $m \geq 1$
satisfying the conditions 
\begin{eqnarray*}
&&0\le l_1\le \frac{n}{2},\\
&&1\le l_i\le \frac{n}{2}\qquad \hbox{for }2\le i\le m
\end{eqnarray*}
and
$$
l_1+\cdots +l_m\le n.
$$
To estimate this term, we choose real numbers $q_1, \ldots, q_m \in [2,\infty[$ such that 
$$
l_1\leq\frac{n}{q_1},\qquad \frac{n}{q_i} <l_i\enspace\hbox{ for } \enspace 2 \leq
i \leq m$$
and $$
 \frac{1}{2}  \leq  \frac{1}{q_1} +\cdots +\frac{1}{q_m} < 1.
$$
Moreover, we define real numbers $\rho_2, \ldots, \rho_m$ by 
$$ \rho_1 = \frac{l_1 - \frac{n}{q_1} + \frac{n}{2}}{\frac{n}{2}} \in [0,1] $$
and 
$$ \rho_i = \frac{l_i - \frac{n}{q_i}}{\frac{n}{2}} \in \hspace{1.2mm} ]0,1[\quad\hbox{
for }2 \leq i \leq m.$$ Then we have 
 $ \rho_1 +\cdots +\rho_m \leq 2$; 
hence 
 $ \rho_2 +\cdots +\rho_m \leq 2 - \rho_1.$ 
Since $w$ is bounded in $H^{\frac{n}{2}}$, this implies  
\begin{eqnarray*} 
&&-\int_M 2 \, \big ( (-\Delta_0)^{\frac{n}{4}} v \big )^2 \, dV_0 
+ C \sum_{l_1, \ldots, l_m} \int_M |\nabla_0^{l_1} v| \cdot 
|\nabla_0^{l_2} w| \cdots |\nabla_0^{l_m} w| \, e^{\alpha w} \, dV_0 
\\ &\leq& -\|v\|_{H^{\frac{n}{2}}}^2 + C 
\sum_{l_1, \ldots, l_m} \|\nabla_0^{l_1} v\|_{L^{q_1}} \cdot 
\|\nabla_0^{l_2} w\|_{L^{q_2}} 
\cdots \|\nabla_0^{l_m} w\|_{L^{q_m}} \\ 
&\leq& -\|v\|_{H^{\frac{n}{2}}}^2 + C 
\sum_{l_1, \ldots, l_m} \|v\|_{H^{l_1 - \frac{n}{q_1} + \frac{n}{2}}} 
\cdot \|w\|_{H^{l_2 - \frac{n}{q_2} + \frac{n}{2}}} 
\cdots \|w\|_{H^{l_m - \frac{n}{q_m} + \frac{n}{2}}} \\ 
&\leq& -\|v\|_{H^{\frac{n}{2}}}^2 + C 
\sum_{l_1, \ldots, l_m} \|v\|_{L^2}^{1-\rho_1} \, 
\|v\|_{H^{\frac{n}{2}}}^{\rho_1} \, \|w\|_{H^{\frac{n}{2}}}^{(1-\rho_2)
+\cdots +(1-\rho_m)} \, 
\|w\|_{H^n}^{\rho_2+\cdots +\rho_m} \\ 
&\leq& -\|v\|_{H^{\frac{n}{2}}}^2 + C 
\sum_{l_1, \ldots, l_m} \|v\|_{L^2}^{1-\rho_1} \, \|v\|_{H^{\frac{n}{2}}}^{\rho_1} \, \|w\|_{H^n}^{\rho_2+\cdots +\rho_m} \\ 
&\leq &C 
\sum_{l_1, \ldots, l_m} \|v\|_{L^2}^{\frac{2-2\rho_1}{2-\rho_1}} \, 
\|w\|_{H^n}^{\frac{2(\rho_2+\cdots +\rho_m)}{2-\rho_1}} \\ 
&\leq& C \, (\|v\|_{L^2}^2 + 1) \, (\|w\|_{H^n}^2 + 1). 
\end{eqnarray*} 
Thus, we conclude that 
$$ \frac{d}{dt} \bigg ( \int_M (P_0 w)^2 \, dV_0 \bigg ) 
\leq C \, (\|v\|_{L^2}^2 + 1) \, (\|w\|_{H^n}^2 + 1).$$  
From the positivity of $P_0$ it follows that 
$$ \|w\|_{H^n}^2 \leq C \int_M (P_0 w)^2 \, dV_0.$$  
Moreover, the function $v$ satisfies 
$$ \|v\|_{L^2}^2 = \int_M \frac{1}{4} \, e^{nw} \, \Big ( 
Q - \frac{\overline{Q} \, f}{\overline{f}} \Big )^2 \, dV_0 = \int_M
\frac{1}{4} \, \Big ( Q - \frac{\overline{Q} \, f}{\overline{f}} \Big )^2 \, dV.$$  
Therefore, we obtain 
$$
\frac{d}{dt} \bigg ( \int_M (P_0 w)^2 \, dV_0 + 1 \bigg )
 \leq C \, \bigg ( \int_M (P_0 w)^2 \, dV_0 +
1 \bigg ) \, \bigg ( \int_M
\Big ( Q - \frac{\overline{Q} \, f}{\overline{f}} \Big )^2 \, dV + 1 \bigg ). 
$$ Since 
$$ \int_0^T \int_M \Big ( Q - \frac{\overline{Q} \, f}{\overline{f}} \Big )^2 \,
dV \, dt \leq C,$$  we deduce that 
$$ \int_M (P_0 w)^2 \, dV_0 \leq C\quad\hbox{ for all }\quad 0 \leq t \leq T.$$ 
This implies $$ \|w\|_{H^n} \leq C$$  for all $0 \leq t \leq T$. 
Using the Sobolev inequality, we obtain 
$$ |w| \leq C\quad\hbox{ for all }\quad 0 \leq t \leq T.$$

\section{Boundedness of $w$ in $H^{2k}$ for $0 \leq t \leq T$}

We now establish bounds for the higher order derivatives:
\begin{eqnarray*} 
\frac{d}{dt} \bigg ( \int_M |(-\Delta_0)^k w|^2 \, dV_0 \bigg ) 
&\leq& -\int_M e^{-nw} \, |(-\Delta_0)^{k+\frac{n}{4}} w|^2 \, dV_0 
\\ &&+ C \sum_{k_1, \ldots, k_m} \int_M |\nabla_0^{k_1} w| \cdots |\nabla_0^{k_m}
w| \, dV_0\,; 
\end{eqnarray*}
hence 
\begin{eqnarray*} 
\frac{d}{dt} \bigg ( \int_M |(-\Delta_0)^k w|^2 \, dV_0 \bigg ) 
&\leq &-\frac{1}{C} \int_M |(-\Delta_0)^{k+\frac{n}{4}} w|^2 \, dV_0 
\\ &&+ C \sum_{k_1, \ldots, k_m} \int_M |\nabla_0^{k_1} w| \cdots |\nabla_0^{k_m}
w| \, dV_0. 
\end{eqnarray*} 
Here, the sum is taken over all $m$-tuples $k_1, \ldots, k_m$, with $m \geq 3$,
which satisfy the conditions 
$$ 1 \leq k_i \leq 2k + \frac{n}{2}\quad\hbox{ and }\quad  k_1 + \cdots + k_m \leq 4k + n.$$  
We now choose real numbers $p_1, \ldots, p_m \in [2,\infty[$ such that 
$$ k_i \leq 2k + \frac{n}{p_i}\hensp{and}\frac{1}{p_1} + \cdots + \frac{1}{p_m} = 1.$$  
Moreover, we define real numbers $\theta_1, \ldots, \theta_m$ by 
$$ \theta_i ={\rm max} \left\{ \frac{k_i-\frac{n}{p_i}-\frac{n}{2}}{2k-\frac{n}{2}} ,0\right\}.$$  
Since $m \geq 3$, we can choose $p_1,\ldots ,p_m\in [2,\infty[$ such that  
$$ \theta_1 + \cdots + \theta_m  < 2.$$  
From this it follows that 
\begin{eqnarray*} 
\frac{d}{dt} \|w\|_{H^{2k}}^2 
&\leq& -\frac{1}{C} \, \|w\|_{H^{2k+\frac{n}{2}}}^2 
+ C \sum_{k_1, \ldots, k_m} \|\nabla_0^{k_1} w\|_{L^{p_1}} \cdots 
\|\nabla_0^{k_m} w\|_{L^{p_m}} \\ 
&\leq& -\frac{1}{C} \, \|w\|_{H^{2k+\frac{n}{2}}}^2 
+ C \sum_{k_1, \ldots, k_m} \|w\|_{H^{k_1-\frac{n}{p_1}+\frac{n}{2}}} \cdots 
\|w\|_{H^{k_m-\frac{n}{p_m}+\frac{n}{2}}} \\ 
&\leq& -\frac{1}{C} \, \|w\|_{H^{2k+\frac{n}{2}}}^2 
+ C \sum_{k_1, \ldots, k_m} \|w\|_{H^n}^{(1-\theta_1) + \cdots + (1-\theta_m)} 
\, \|w\|_{H^{2k+\frac{n}{2}}}^{\theta_1 + \cdots + \theta_m} \\ 
&\leq& -\frac{1}{C} \, \|w\|_{H^{2k+\frac{n}{2}}}^2 
+ C \sum_{k_1, \ldots, k_m} \|w\|_{H^{2k+\frac{n}{2}}}^{\theta_1 + \cdots + \theta_m} \\ 
&\leq &-\frac{1}{C} \, \|w\|_{H^{2k+\frac{n}{2}}}^2 + C \\ 
&\leq &-\frac{1}{C} \, \|w\|_{H^{2k}}^2 + C 
\end{eqnarray*} 
for all $0 \leq t \leq T$. Thus, we conclude that $$ \|w\|_{H^{2k}} \leq C\hensp{ 
for all }0 \leq t \leq T.
$$
Therefore, the evolution equation has a 
solution which is defined for all time.

\section{Convergence}

For the sake of brevity, we put 
$$ y(t) = \int_M \Big ( Q - \frac{\overline{Q} \, f}{\overline{f}} \Big )^2 \, dV $$  
and we show that $$ y(t) \to 0\hensp{for}t \to \infty.
$$
 Let $\varepsilon$ be an arbitrary positive
number. We choose $t_0 \geq 0$ such that  $ y(t_0) \leq \varepsilon.$   We claim 
that  $ y(t) \leq 3\varepsilon$   for all $t \geq t_0$. Otherwise, we define 
$$ t_1 =
\inf \{t \geq t_0: y(t) \geq 3\varepsilon\}.$$  This implies 
$  y(t) \leq 3\varepsilon$  for all $t_0 \leq t \leq t_1$. 
From this it follows that $$ \int_M e^{-nw} \, (Q_0 + P_0 w)^2 \, dV_0 \leq C$$  for
all $t_0 \leq t \leq t_1$. Moreover, it follows from   results 
in Section 3 that 
$$ \int_M e^{3nw} \, dV_0 \leq C\hensp{for all} t_0 \leq t \leq t_1.$$
Using H\"older's inequality, we obtain 
$$
 \int_M |Q_0 + P_0 w|^{\frac{3}{2}} \, dV_0  \leq \bigg ( \int_M e^{-nw} \, (Q_0 + P_0 w)^2 \, dV_0 \bigg
)^{\frac{3}{4}} \, 
\bigg ( \int_M e^{3nw} \, dV_0 \bigg )^{\frac{1}{4}}. $$
From this it follows that 
$$ \int_M |P_0 w|^{\frac{3}{2}} \, dV_0 \leq C\hensp{for all } t_0 \leq t \leq t_1.$$ Using the Sobolev inequality, we
obtain 
$$ |w| \leq C\hensp{for all}t_0 \leq t \leq t_1.
$$

We have shown in Section 2 that the function $Q - \frac{\overline{Q} \,
f}{\overline{f}}$ satisfies the evolution equation 
\begin{eqnarray*} 
\frac{\partial}{\partial t} \Big ( Q - \frac{\overline{Q} \, f}{\overline{f}} 
\Big ) &= &-\frac{1}{2} \, P \Big ( Q - \frac{\overline{Q} \, f}{\overline{f}} 
\Big ) + \frac{n}{2} \, Q \Big ( Q - \frac{\overline{Q} \, f}{\overline{f}} \Big ) 
\\ &&- \frac{n}{2} \, \frac{\overline{Q} \, f}{\overline{f}} \int_M \frac{f}{\overline{f}} \, 
\Big ( Q - \frac{\overline{Q} \, f}{\overline{f}} \Big ) \, 
dV, \end{eqnarray*} 
where $P$ denotes the Paneitz
operator with respect to the metric $g$. From this it follows that 
\begin{eqnarray*} 
 \frac{d}{dt} \bigg ( \int_M \Big ( 
Q - \frac{\overline{Q} \, f}{\overline{f}} \Big )^2 \, dV \bigg ) 
&= &-\int_M \Big ( Q - \frac{\overline{Q} \, f}{\overline{f}} \Big
) \, P \Big ( Q - \frac{\overline{Q} \, f}{\overline{f}} \Big ) \, dV 
\\ &&+ \int_M \frac{n}{2} \, \Big ( Q - \frac{\overline{Q} \, f}{\overline{f}} \Big )^3 \, dV . 
\\ &&+ \int_M n \, \frac{\overline{Q} \, f}{\overline{f}} \, 
\Big ( Q - \frac{\overline{Q} \, f}{\overline{f}} \Big )^2 \, dV 
\\ &&- n \, \overline{Q} \, \bigg ( \int_M \frac{f}{\overline{f}} \, 
\Big ( Q - \frac{\overline{Q} \, f}{\overline{f}} \Big ) \, dV 
\bigg )^2. 
\end{eqnarray*} 
Using the Gagliardo-Nirenberg inequality, we can bound 
$$ \Big \| Q - \frac{\overline{Q} \, f}{\overline{f}} \Big \|_{L^3} 
\leq C \, \Big \| Q - \frac{\overline{Q} \, f}{\overline{f}} 
\Big \|_{L^2}^{\frac{2}{3}} 
\, \Big \| Q - \frac{\overline{Q} \,
f}{\overline{f}} \Big \|_{H^{\frac{n}{2}}}^{\frac{1}{3}},$$  
where the norms are taken with respect to the background metric $g_0$. This implies 
\begin{eqnarray*} 
&&\int_M 
\Big ( Q - \frac{\overline{Q} \, f}{\overline{f}} \Big )^3 \, dV_0 \\ 
&&\qquad \leq C \, \bigg ( \int_M \Big ( Q - \frac{\overline{Q} \, f}{\overline
{f}} \Big )^2 \, dV_0 \bigg ) \, 
\bigg ( \int_M \Big ( Q - \frac{\overline{Q} \, f}{\overline{f}} \Big ) \, 
P_0 \Big ( Q - \frac{\overline{Q} \, f}{\overline{f}} \Big ) \, dV_0 \bigg
)^{\frac{1}{2}}. \end{eqnarray*} 
Since $w$ is uniformly bounded for $t_0 \leq t \leq t_1$, we obtain 
\begin{eqnarray*} 
&&\int_M 
\Big ( Q - \frac{\overline{Q} \, f}{\overline{f}} \Big )^3 \, dV \\ 
&&\qquad \leq C \, \bigg ( \int_M \Big ( Q - \frac{\overline{Q} \, f}{\overline
{f}} \Big )^2 \, dV \bigg ) \, 
\bigg ( \int_M \Big ( Q - \frac{\overline{Q} \, f}{\overline{f}} \Big ) \, 
P \Big ( Q - \frac{\overline{Q} \, f}{\overline{f}} \Big ) \, dV \bigg
)^{\frac{1}{2}}. \end{eqnarray*} 
Thus, we conclude that 
\begin{eqnarray*} 
 \frac{d}{dt} \bigg ( \int_M \Big ( Q - \frac{\overline
{Q} \, f}{\overline{f}} \Big )^2 \, dV \bigg )  &\leq& C \, 
\bigg ( \int_M \Big ( Q - \frac{\overline{Q} \, f}{\overline{f}} \Big )^2 
\, dV \bigg )^2\\
&& + C \, \bigg ( \int_M 
\Big ( Q - \frac{\overline{Q} \, f}{\overline{f}} \Big )^2 \, dV \bigg ); 
\end{eqnarray*} 
hence $$ \frac{d}{dt} y(t) \leq C \, y(t)^2 + C \, y(t).$$  
Therefore, we obtain 
$$ 2\varepsilon \leq y(t_1) - y(t_0) \leq C \int_{t_0}^{t_1} y(t) \, dt.$$  
If we choose $t_0$ large enough, then we have 
$$ C \int_{t_0}^\infty y(t) \, dt \leq \varepsilon.$$  Hence, we obtain $ 2\varepsilon
\leq \varepsilon $ which is a contradiction. 
Thus, we conclude that $$ y(t) \to 0\hensp{ for }t \to \infty.
$$
From this it follows that 
$$ |w| \leq C\hensp{for all} t \geq 0.
$$
 Moreover, we have 
$$ \int_M e^{-nw} \, (Q_0 + P_0 w)^2 \, dV_0 \leq C$$  for all $t \geq 0$. 
From this it follows that 
$$ \int_M (Q_0 + P_0 w)^2 \, dV_0 \leq C;$$  hence $$ \|w\|_{H^n} \leq C$$  for all $t 
\geq 0$. Arguing as above, we obtain $$ \|w\|_{H^k} \leq C\hensp{for all} t \geq 0.$$

It remains to show that the flow converges to a metric 
satisfying $$ \frac{Q}{f} = \frac{\overline{Q}}{\overline{f}}.$$  
The evolution equation $$ \frac{\partial}{\partial t} g = 
- \Big ( Q - \frac{\overline{Q} \, f}{\overline{f}} \Big ) \, g$$  
is the gradient flow for the functional $$ E_f[w] = \int_M \frac{n}{2} \, w \, P_0 w \, dV_0 + \int_M n \, Q_0 \, w \, dV_0 - 
\int_M Q_0 \, dV_0 \,
\log \bigg ( \int_M e^{nw} \, f \, dV_0 \bigg ).$$  Since the functional 
$E_f[w]$ is real analytic, the assertion follows from a general result 
of L.\ Simon \cite{Si}.

\section{The case $M = {\bf RP}^n$}

In this section, we consider the special case $M = {\bf RP}^n$. 
We normalize the metric such that 
the volume of $M$ is equal to $\frac{1}{2} \, \omega_n$ and 
the mean value of the function $Q$ is equal to $(n-1)!$. 
By Theorem \ref{convergence.1}, the flow converges to a limit metric $g$ satisfying 
$$ \frac{Q}{f} = \frac{(n-1)!}{\overline{f}}.$$  
In particular, for every positive function $f$ on ${\bf RP}^n$, there exists a metric $g$ on 
${\bf RP}^n$ such that 
$$ \frac{Q}{f} = \frac{(n-1)!}{\overline{f}}.$$  
We now consider the case $f = 1$. In this case, the limit metric $g$ satisfies 
$Q = (n-1)!$. It follows from a result of S.-Y.\ A.\  Chang and P.\ 
Yang \cite{CY2}
(see also C.\ S.\  Lin's paper \cite{Li}) 
that the limit metric is the standard metric on~${\bf RP}^n$. 

We claim that the flow converges exponentially. To show this, we denote by $g_0$ 
the standard metric on ${\bf RP}^n$. Then the conformal factor satisfies the
evolution equation $$ \frac{\partial}{\partial t} w = -\frac{1}{2} \, e^{-nw} \, 
P_0 w + \frac{1}{2} \, (n-1)! \, (1 - e^{-nw}).$$  
Linearizing this equation, we obtain $$ \frac{\partial}{\partial t} w = -\frac{1}{2} 
\, P_0 w + \frac{1}{2} \, n! \, w.$$  
The Paneitz operator on ${\bf RP}^n$ is given by 
$$ P_0 = \prod_{k=1}^{\frac{n}{2}} (-\Delta_0 + (k-1)(n-k)).$$  
The first eigenvalue of the Laplace operator $-\Delta_0$ on ${\bf RP}^n$ is
strictly greater than $n$. Hence, the first eigenvalue of the Paneitz operator $P_0$ is strictly greater than $n!$. 
Therefore, the first eigenvalue of 
the linearized operator 
is strictly less than $0$. Thus, we conclude that the flow 
converges exponentially to the standard metric on ${\bf RP}^n$.  

\section{A compactness result for conformal metrics on $S^n$}

In this section, we give a proof for Proposition \ref{compactness.result}. 
Let $g_k = e^{2w_k} \, g_0$ be a sequence of conformal metrics on $S^n$ with fixed
volume such that $$ \int_{S^n} Q_k^2 \, dV_k \leq C.$$  Since 
$$ Q_k = e^{-nw_k} \, (Q_0 + P_0w_k),$$  we obtain 
$$ \int_{S^n} e^{-nw_k} \, (Q_0 + P_0 w_k)^2 \, dV_0 \leq C.$$  
Moreover, we have $$ \int_{S^n} |Q_k| \, dV_k \leq C.$$  
Hence $$ \int_{S^n} |P_0 w_k| \, dV_0 \leq C.$$  
Finally, we have 
$$ \lim_{r \to 0} \lim_{k \to \infty} \int_{B_r(x)} |Q_k| \, dV_k < \frac{1}{2}
\, (n-1)! \, \omega_n.$$  
This implies 
$$ \lim_{r \to 0} \lim_{k \to \infty} \int_{B_r(x)} |P_0 w_k| \, dV_0 <
\frac{1}{2} \, (n-1)! \, \omega_n.$$  
Choosing $r$ sufficiently small, we obtain 
$$ \lim_{k \to \infty} \int_{B_r(x)} |P_0 w_k| \, dV_0 < \frac{1}{2} \, (n-1)! \,
\omega_n.$$  
Let $$ I_k = \int_{B_r(x)} |P_0 w_k| \, dV_0.$$  
We now use the formula 
$$ w_k(y) - \overline{w}_k = 
\int_{S^n} P_0 w_k(z) \, K(y,z) \, dV_0(z).$$  
This implies 
$$ np (w_k(y) - \overline{w}_k) \leq \int_{B_r(x)} np \, |P_0 w_k(z)| \, 
|K(y,z)| \, dV_0(z) + C$$  
for all $y \in B_{\frac{r}{2}}(x)$. Using Jensen's inequality, we obtain 
$$ e^{np (w_k(y) - \overline{w}_k)} \leq \frac{C}{I_k} \int_{B_r(x)} |P_0 
w_k(z)| \, e^{npI_k \, |K(y,z)|} \, dV_0(z)$$  
for all $y \in B_{\frac{r}{2}}(x)$. Since 
$$ \lim_{k \to \infty} I_k < \frac{1}{2} \, (n-1)! \, \omega_n,$$  we can find a real number $p > 1$
such that $$ \lim_{k \to \infty} pI_k < \frac{1}{2} \, (n-1)! \, \omega_n.$$  
We now use an asymptotic formula of the
function $K(y,z)$ for $|y - z| \to 0$. To derive this formula, we use the
identity $$ (-\Delta)^{\frac{n}{2}} \log |y - z| = -2^{n-2} \, \big ( \big ( \frac
{n-2}{2} \big )! \big )^2 \, \omega_{n-1} \, \delta(y - z).$$  
This implies 
$$ (-\Delta)^{\frac{n}{2}} \log |y - z| = -\frac{1}{2} \, (n-1)! \, \omega_n 
\, \delta(y - z).$$  Therefore, we obtain 
$$ \frac{1}{2} \, (n-1)! \, \omega_n \, K(y,z) \sim -\log |y - z|;$$  
hence $$ e^{\frac{1}{2} \, (n-1)! \, \omega_n \, |K(y,z)|} \sim 
\frac{1}{|y-z|}.$$  
From this it follows that $$ \int_{S^n} e^{npI_k \, |K(y,z)|} \, dV_0(y) 
\leq C.$$  Since 
$$ \frac{1}{I_k} \int_{B_r(x)} |P_0 w_k| \, dV_0 = 1,$$  we conclude that 
$$ \int_{B_{\frac{r}{2}}(x)} 
e^{np (w_k(y) - \overline{w}_k)} \, dV_0(y) \leq C.$$  
Covering $S^n$ with finitely many balls $B_{\frac{r}{2}}(x)$, we
obtain $$ \int_{S^n} e^{np(w_k - \overline{w}_k)} \, dV_0 \leq C$$  for some 
$p > 1$. In particular, we have $$ \int_{S^n} e^{n(w_k - \overline{w}_k)} \, dV_0 \leq
C.$$  Since $ \int_{S^n} e^{nw_k} \, dV_0 = 1,$   we conclude that
 $ e^{-n\overline{w}_k} \leq C;$ 
hence 
 $ -C \leq \overline{w}_k \leq C.$ 
This implies 
$$ \int_{S^n} e^{npw_k} \, dV_0 \leq C.$$  
By H\"older's inequality,  
$$\int_{S^n} |Q_0 + P_0 w_k|^{\frac{2p}{p+1}} \, dV_0 \leq \bigg ( \int_{S^n} e^{-nw_k} \, (Q_0 + P_0 w_k)^2 \, dV_0
\bigg )^{\frac{p}{p+1}} \, 
\bigg ( \int_{S^n} e^{npw_k} \, dV_0 \bigg )^{\frac{1}{p+1}}. $$
From this it follows that 
$$ \int_{S^n} |P_0 w_k|^{\frac{2p}{p+1}} \, dV_0 \leq C.$$  
Using the Sobolev inequality, we obtain 
$ |w_k| \leq C.$  Thus, we conclude that $$ \int_{S^n} |P_0 w_k|^2 \, dV_0
\leq C.$$  Therefore, the sequence $w_k$ is uniformly bounded in $H^n$.

 \end{document}